\documentclass{scrartcl}
\usepackage{amssymb}
\usepackage{amsmath}
\usepackage{amsfonts}
\usepackage{latexsym}
\usepackage{color,graphicx}  
\usepackage[bottom]{footmisc}
\usepackage{cite}
\usepackage{fancyvrb}  
\usepackage{times}
\usepackage{mathptmx}
\usepackage{hyperref}

\newcommand{\bbbr}{\mathbb{R}}
\newcommand{\Idx}{\mathcal{I}}
\newcommand{\supp}{\mathop{\operatorname{supp}}\nolimits}
\newcommand{\diam}{\mathop{\operatorname{diam}}\nolimits}
\newcommand{\dist}{\mathop{\operatorname{dist}}\nolimits}

\begin{document}

\title{GCA-$\mathcal{H}^2$ matrix compression for electrostatic
  simulations}
\author{Steffen B\"orm \and Sven Christophersen}

\maketitle

\begin{abstract}
We consider a compression method for boundary element matrices arising
in the context of the computation of electrostatic fields.
Green cross approximation combines an analytic approximation of the
kernel function based on Green's representation formula and quadrature
with an algebraic cross approximation scheme in order to obtain both the
robustness of analytic methods and the efficiency of algebraic ones.
One particularly attractive property of the new method is that it is
well-suited for acceleration via general-purpose graphics processors
(GPUs).
\end{abstract}

\section{Introduction}
Boundary integral formulations are particularly useful when dealing
with electrostatic exterior domain problems:
we only have to construct a mesh for the boundary of the domain, and
once an integral equation on this boundary has been solved, we can
directly evaluate the electrostatic field in all points of the infinite
domain by computing a surface integral.

Standard formulations typically lead to equations of the form
\begin{equation*}
  \int_{\partial\Omega} g(x,y) u(y) \,dy
  = \lambda u(x)
    + \int_{\partial\Omega} \frac{\partial g}{\partial n_y}(x,y) v(y) \,dy
\end{equation*}
for all $x\in\partial\Omega$, where $\Omega\subseteq\bbbr^3$ is a domain,
$\lambda\in\bbbr$, $u$ and $v$ are scalar functions on the boundary
$\partial\Omega$, and
\begin{equation*}
  g(x,y) = \frac{1}{4\pi \|x-y\|}
\end{equation*}
is the fundamental solution of Laplace's equation.

Discretization by Galerkin's method with basis functions
$(\varphi_i)_{i\in\Idx}$ leads to a matrix $G\in\bbbr^{\Idx\times\Idx}$ given
by
\begin{equation}\label{eq:G_def}
  g_{ij} = \int_{\partial\Omega} \varphi_i(x)
  \int_{\partial\Omega} g(x,y) \varphi_j(y) \,dy \,dx
\end{equation}
for all $i,j\in\Idx$, and all of these coefficients are typically non-zero.

Working directly with the matrix $G$ is unattractive, since for
$n:=\#\Idx$ basis functions, we would have to store $n^2$ coefficients
and quickly run out of memory.

This problem can be fixed by taking advantage of the properties of the
kernel function $g$:
analytic approximation schemes like the fast multipole method
\cite{RO85,GRRO97}, Taylor expansion \cite{HANO89}, or interpolation
\cite{BOGR02,BOLOME02} replace $g$ in suitable subdomains of the boundary
$\partial\Omega$ by a short sum
\begin{equation*}
  g(x,y) \approx \sum_{\nu=1}^k a_\nu(x) b_\nu(y)
\end{equation*}
that leads to a low-rank approximation of corresponding submatrices of $G$,
while algebraic schemes like the adaptive cross approximation (ACA)
\cite{TY96,BE00a,BERJ01} or rank-revealing factorizations
\cite{CHIP94} directly construct low-rank approximations based
on the matrix entries.

Hybrid approximation schemes like generalized fast multipole methods
\cite{GIRO02,BIYIZO04} or hybrid cross approximation (HCA) \cite{BOGR04}
combine the concepts of analytic and algebraic approximation in order to
obtain the near-optimal compression rates of algebraic methods while
preserving the stability and robustness of analytic techniques.

Our algorithm falls into the third category:
\emph{Green cross approximation} (GCA) combines an analytic approximation
based on Green's representation formula with adaptive cross approximation
(ACA) to obtain low-rank approximations of submatrices.
In order to improve the efficiency, we employ GCA in a recursive fashion
that allows us to significantly reduce the storage requirements without
losing the method's fast convergence.

While all of these technique allow us to handle the boundary integral
equation more or less efficiently, analytic methods and some of the
hybrid methods can also be used to speed up the subsequent evaluation
of the electrostatic field in arbitrary points of $\Omega$.

\section{Green quadrature}
In order to find a data-sparse approximation of $G$, we consider
a domain $\tau\subseteq\bbbr^3$ and a superset $\omega\subseteq\bbbr^3$
such that the distance from $\tau$ to the boundary $\partial\omega$
of $\omega$ is non-zero.
For any $y\in\bbbr^3\setminus\overline{\omega}$, the function
$x \mapsto g(\cdot,y)$ is harmonic in $\omega$, so we can
apply Green's representation formula (also known as Green's third identity) to obtain
\begin{equation*}
  g(x,y)
  = \int_{\partial\omega}
       g(x,z) \frac{\partial g}{\partial n_z}(z,y)
     - \frac{\partial g}{\partial n_z}(x,z) g(z,y) \,dz
\end{equation*}
for all $x\in\tau$ and $y\in\bbbr^3\setminus\overline{\omega}$.
If the distances between $\partial\omega$ and $\tau$ and between
$\partial\omega$ and $y$ are sufficiently large, the integrand is
smooth, and we can approximate the integral by a quadrature rule
to find
\begin{align}
  g(x,y)
  \approx \sum_{\nu=1}^k
      &w_\nu g(x,z_\nu) \frac{\partial g}{\partial n_z}(z_\nu,y)\notag\\
    - &w_\nu \frac{\partial g}{\partial n_z}(x,z_\nu) g(z_\nu,y)
    \label{eq:quadrature_apx}
\end{align}
with weights $w_\nu$ and quadrature points $z_\nu$, and in this
approximation the variables $x$ and $y$ are separated.

This gives rise to a first low-rank approximation of $G$:
given subsets $\widehat{\tau},\widehat{\sigma}\subseteq\Idx$ of the
index set, we can introduce axis-parallel boxes
\begin{align*}
  \tau &\supseteq \bigcup_{i\in\widehat{\tau}} \supp \varphi_i, &
  \sigma &\supseteq \bigcup_{j\in\widehat{\sigma}} \supp \varphi_j
\end{align*}
containing the supports of the corresponding basis functions,
and if these boxes are well-separated, we can find a superset
$\omega$ of $\tau$ such that its boundary $\partial\omega$
is sufficiently far from both $\tau$ and $\sigma$.
Replacing $g$ in the definition (\ref{eq:G_def}) of the Galerkin
matrix by the quadrature-based approximation
leads to a factorized approximation
\begin{equation*}
  G|_{\widehat{\tau}\times\widehat{\sigma}}
  \approx A_{\tau\sigma} B_{\tau\sigma}^*,
\end{equation*}
with $A_{\tau\sigma}\in\bbbr^{\widehat{\tau}\times 2k}$ and
$B_{\tau\sigma}\in\bbbr^{\widehat{\sigma}\times 2k}$, so the
rank of the approximation is bounded by $2k$.

The matrix coefficients are given by
\begin{align*}
  a_{\tau\sigma,i\nu}
  &= \sqrt{w_\nu} \int_{\partial\Omega} g(x,z_\nu) \varphi_i(x) \,dx,\\
  a_{\tau\sigma,i(\nu+k)}
  &= -d_\tau \sqrt{w_\nu} \int_{\partial\Omega} \frac{\partial g}{\partial n_z}(x,z_\nu)
  \varphi_i(x) \,dx,\\
  b_{\tau\sigma,j\nu}
  &= \sqrt{w_\nu} \int_{\partial\Omega} \frac{\partial g}{\partial n_z}(z_\nu,y)
  \varphi_j(y) \,dy,\\
  b_{\tau\sigma,j(\nu+k)}
  &= \frac{\sqrt{w_\nu}}{d_\tau} \int_{\partial\Omega} g(z_\nu,y) \varphi_j(y) \,dy,
\end{align*}
where the scaling parameter $d_\tau = \diam(\tau)$ serves to balance the
different scaling behaviour of the kernel function and its normal derivative.

%
%
\begin{figure}
  \begin{center}
  \includegraphics[width=0.8\textwidth]{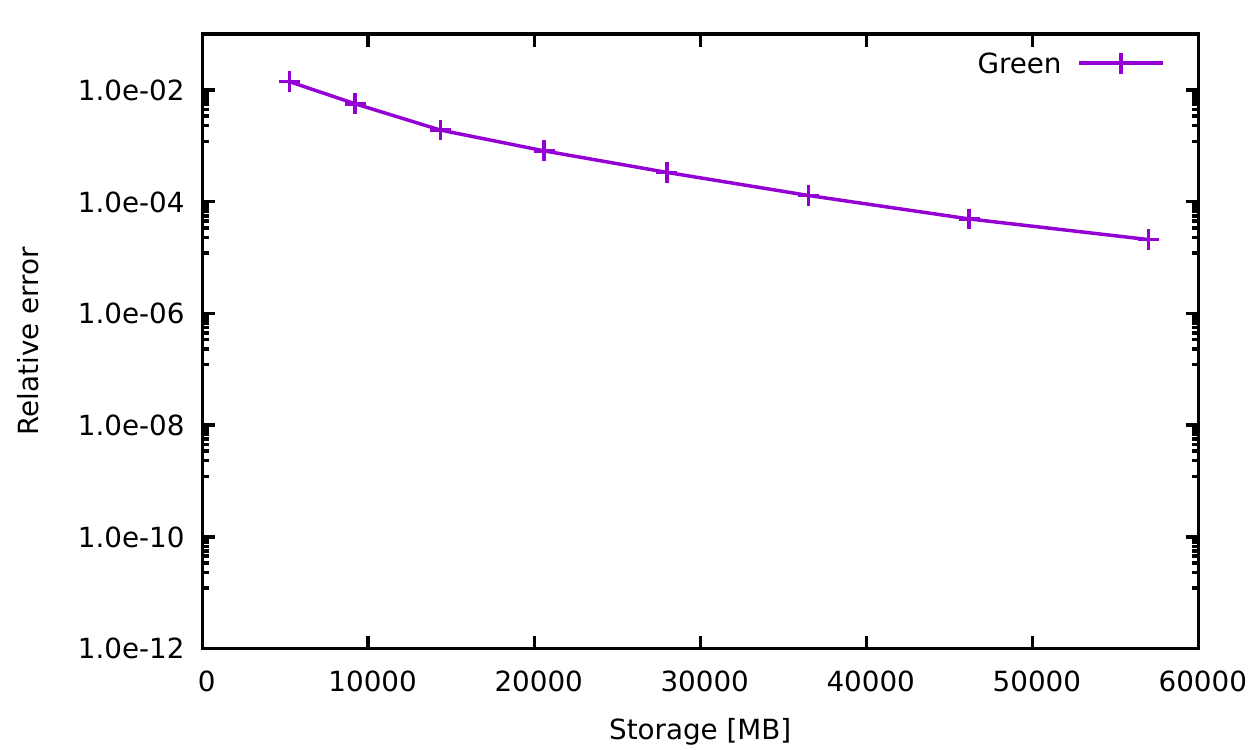}
  \caption{Relative error of the Green quadrature approximation compared to
    the storage requirements}
  \label{fi:evm_green}
  \end{center}
\end{figure}

We apply the approximation scheme to a polygonal approximation of the
unit sphere by $n=32\,768$ triangles, choosing piecewise constant basis
functions and the admissibility condition
\begin{equation*}
  \max\{\diam(\tau),\diam(\sigma)\} \leq 2 \eta \dist(\tau,\sigma)
\end{equation*}
with a parameter $\eta\in\bbbr_{>0}$ to determine whether a submatrix
$G|_{\widehat\tau\times\widehat\sigma}$ can be approximated.

Figure~\ref{fi:evm_green} shows the relative spectral error, estimated via
the power iteration, as a function of the storage requirements.
We can see that the convergence is quite disappointing, particularly since
storing the entire matrix as a simple two-dimensional array requires only
$8\,192$ MB of storage.

\section{Green cross approximation}
In order to make the approximation more efficient, we can apply
adaptive cross approximation \cite{BE00a} to derive the algebraic
counterpart of interpolation:
this technique provides us with a small subset
$\tilde\tau\subseteq\widehat\tau$ and a matrix
$V_\tau\in\bbbr^{\widehat\tau\times\tilde\tau}$ such that
\begin{equation*}
  V_\tau A_{\tau\sigma}|_{\tilde\tau\times 2k} \approx A_{\tau\sigma},
\end{equation*}
i.e., we can reconstruct $A_{\tau\sigma}$ using only a few of
its rows.
Since $A_{\tau\sigma}$ is a thin matrix, we can afford to use
reliable pivoting strategies and do not have to rely on heuristics.
We conclude
\begin{equation*}
  V_\tau G|_{\tilde\tau\times\widehat\sigma}
  \approx V_\tau A_{\tau\sigma}|_{\tilde\tau\times 2k} B_{\tau\sigma}^*
  \approx A_{\tau\sigma} B_{\tau\sigma}^*
  \approx G|_{\widehat\tau\times\widehat\sigma},
\end{equation*}
i.e., the algebraic interpolation can also be applied directly to
the original matrix $G$ instead of the low-rank approximation.
This is called a \emph{Green cross approximation} (GCA).

It is important to keep in mind that the matrices $A_{\tau\sigma}$ only
depend on $\tau$, but not on $\sigma$, so the cross approximation
algorithm has to be performed only once for each $\tau$ and both the
set $\tilde\tau$ and the matrix $V_\tau$ do not depend on $\sigma$.

Our modification has two major advantages: on one hand, the
ranks are bounded by both the cardinality of $\widehat\tau$
and the number of quadrature points, so that the approximation
can be far more efficient for small clusters.
On the other hand, we can reach significantly higher accuracies,
since the Green quadrature is only used to choose good
``interpolation points'' $\tilde\tau$, while the final
approximation relies on the entries of the original matrix $G$.

%
%
\begin{figure}
  \begin{center}
  \includegraphics[width=0.8\textwidth]{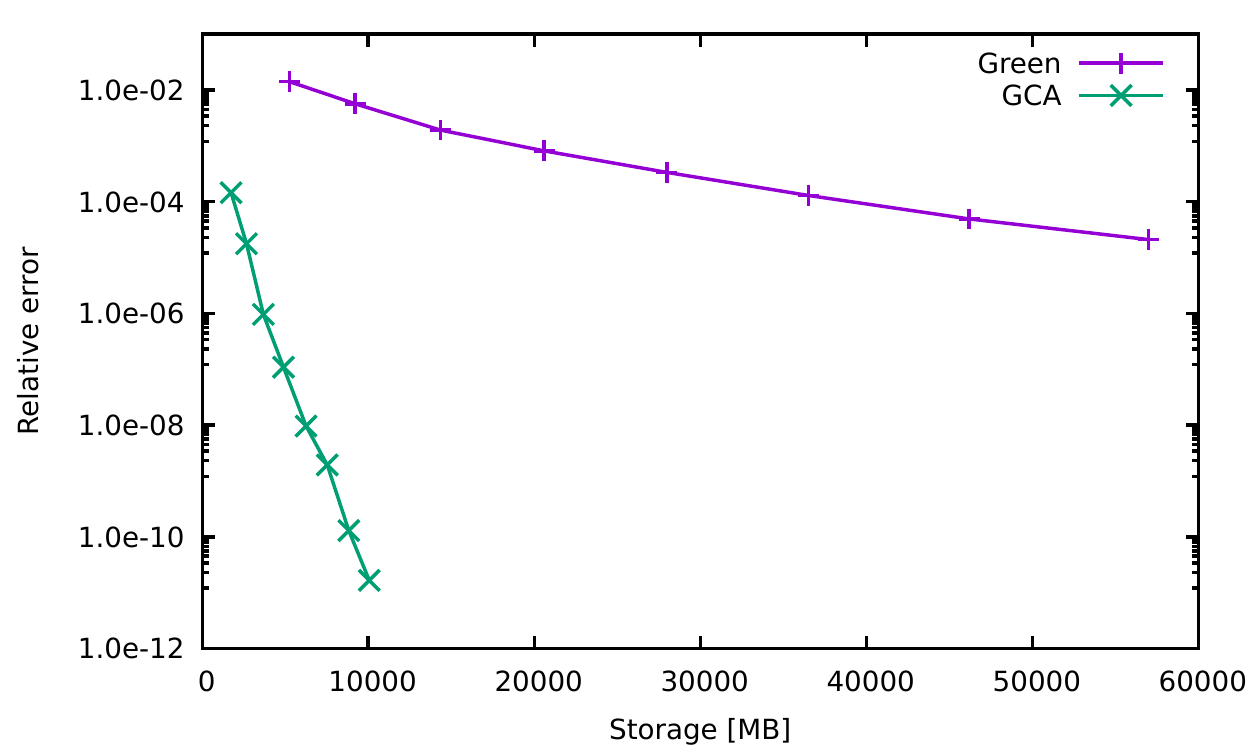}
  \caption{Relative error of the Green quadrature approximation and the
    Green cross approximation (GCA) compared to the storage requirements}
  \label{fi:evm_green_gca}
  \end{center}
\end{figure}

Figure~\ref{fi:evm_green_gca} illustrates that combining cross approximation
with Green quadrature significantly improves the performance: we can reach
fairly high accuracies with moderate storage requirements.

\section{\texorpdfstring{$\mathcal{H}^2$-matrices}{H2-matrices}}
Since Green's formula is symmetric with respect to $\tau$ and $\sigma$,
we can also apply the representation formula to a superset
of $\sigma$ and combine the formular with quadrature and cross
approximation to obtain a subset $\tilde\sigma\subseteq\widehat\sigma$ and
$V_\sigma\in\bbbr^{\widehat\sigma\times\tilde\sigma}$ with
\begin{equation*}
  G|_{\widehat\tau\times\widehat\sigma}
  \approx G|_{\widehat\tau\times\tilde\sigma} V_\sigma^*.
\end{equation*}
Together with the approximation for $\tau$ introduced before, we obtain
the \emph{symmetric} factorization
\begin{equation*}
  G|_{\widehat\tau\times\widehat\sigma}
  \approx V_\tau G|_{\tilde\tau\times\tilde\sigma} V_\sigma^*,
\end{equation*}
and this turns out to be very efficient, since $G|_{\tilde\tau\times\tilde\sigma}$
is usually significantly smaller than $G|_{\widehat\tau\times\widehat\sigma}$.

We can improve the construction further by representing the \emph{basis
matrices} $V_\tau$ and $V_\sigma$ in a hierarchy:
assume that $\widehat\tau$ is subdivided into disjoint subsets
$\widehat\tau_1$ and $\widehat\tau_2$ and that matrices
$V_{\tau_1},V_{\tau_2}$ and subsets
$\tilde\tau_1\subseteq\widehat\tau_1,\tilde\tau_2\subseteq\widehat\tau_2$
have already been constructed.
We let $\tilde\tau_{1,2}:=\tilde\tau_1\cup\tilde\tau_2$ and observe
\begin{align*}
  A_{\tau\sigma}
  &= \begin{pmatrix}
       A_{\tau\sigma}|_{\widehat\tau_1\times 2k}\\
       A_{\tau\sigma}|_{\widehat\tau_2\times 2k}
     \end{pmatrix}
   \approx \begin{pmatrix}
       V_{\tau_1} A_{\tau\sigma}|_{\tilde\tau_1\times 2k}\\
       V_{\tau_2} A_{\tau\sigma}|_{\tilde\tau_2\times 2k}
     \end{pmatrix}
   = \begin{pmatrix}
       V_{\tau_1} & \\
       & V_{\tau_2}
     \end{pmatrix} A_{\tau\sigma}|_{\tilde\tau_{1,2}\times 2k}.
\end{align*}
If we now apply cross approximation to the right factor
$\widehat{A}_{\tau\sigma} := A_{\tau\sigma}|_{\tilde\tau_{1,2}\times 2k}$, we obtain a subset
$\tilde\tau\subseteq\tilde\tau_{1,2}$ and a matrix
$\widehat{V}_\tau\in\bbbr^{\tilde\tau_{1,2}\times\tilde\tau}$ with
\begin{equation*}
  \widehat{V}_\tau \widehat{A}_{\tau\sigma}|_{\tilde\tau\times 2k}
  \approx \widehat{A}_{\tau\sigma}
\end{equation*}
and therefore
\begin{align*}
  A_{\tau\sigma}
  &\approx \begin{pmatrix}
    V_{\tau_1} & \\
    & V_{\tau_2}
  \end{pmatrix} \widehat{A}_{\tau\sigma}
  \approx \begin{pmatrix}
    V_{\tau_1} & \\
    & V_{\tau_2}
  \end{pmatrix} \widehat{V}_\tau A_{\tau\sigma}|_{\tilde\tau\times 2k}
  = V_\tau A_{\tau\sigma}|_{\tilde\tau\times 2k},
\end{align*}
where the basis matrix
\begin{equation*}
  V_\tau := \begin{pmatrix}
    V_{\tau_1} & \\
    & V_{\tau_2}
  \end{pmatrix} \widehat{V}_\tau
\end{equation*}
can be expressed in the form
\begin{align*}
  V_\tau &= \begin{pmatrix}
    V_{\tau_1} E_{\tau_1}\\
    V_{\tau_2} E_{\tau_2}
  \end{pmatrix}, &
  E_{\tau_1} &:= \widehat{V}_\tau|_{\tilde\tau_1\times\tilde\tau}, &
  E_{\tau_2} &:= \widehat{V}_\tau|_{\tilde\tau_2\times\tilde\tau}.
\end{align*}
If we use this factorized representation of the matrices $V_\tau$, we only
have to store $V_\tau$ if $\widehat\tau$ has no subsets, while we use
the substantially smaller \emph{transfer matrices} $E_{\tau_1},E_{\tau_2}$ for
all other index sets.

Since $\tilde\tau_{1,2}$ is usually significantly smaller than
$\widehat\tau$, this construction is faster than the straightforward
GCA approach, and the recursive use of transfer matrices reduces the storage
requirements.
The resulting approximation of $G$ is known as an \emph{$\mathcal{H}^2$-matrix}
\cite{HAKHSA00,BOHA02,BO10}, and it can be proven to have \emph{linear}
complexity with respect to the matrix dimension $n$.

%
%
\begin{figure}
  \begin{center}
  \includegraphics[width=0.8\textwidth]{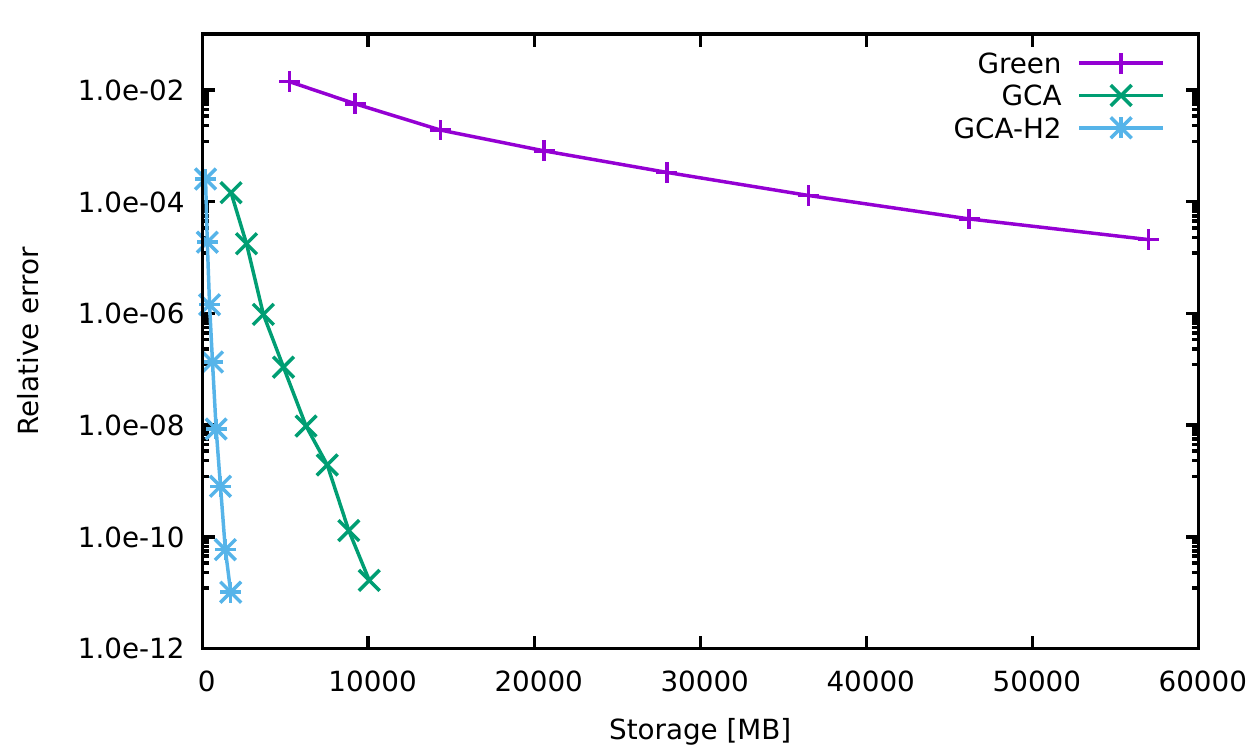}
  \caption{Relative error of the Green quadrature approximation, GCA, and
    GCA-$\mathcal{H}^2$ compared to the storage requirements}
  \label{fi:evm_green_gcah2}
  \end{center}
\end{figure}

\begin{figure}
  \begin{center}
  \includegraphics[width=0.8\textwidth]{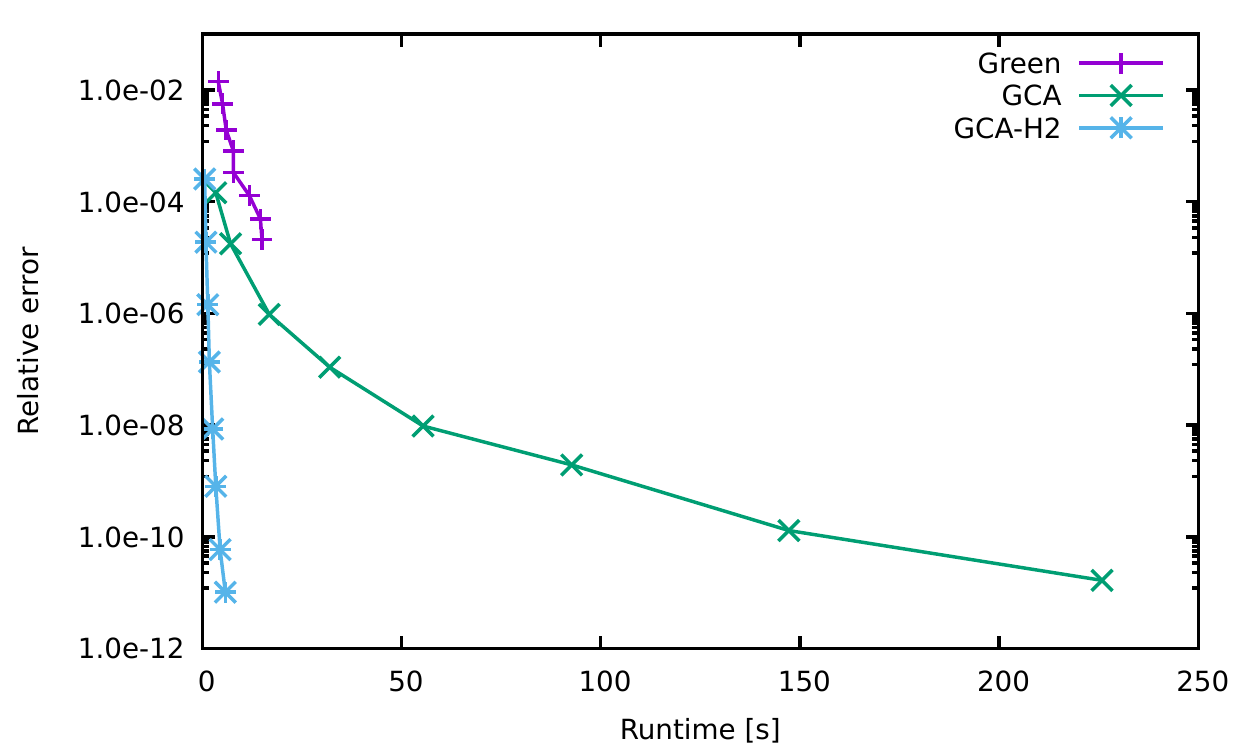}
  \caption{Relative error of the Green quadrature approximation, GCA, and
    GCA-$\mathcal{H}^2$ compared to the setup time}
  \label{fi:evt_green_gcah2}
  \end{center}
\end{figure}

The resulting GCA-$\mathcal{H}^2$-matrix compression algorithm
can be proven to converge exponentially and to have almost optimal
complexity \cite{BOCH14}.
Indeed, Figure~\ref{fi:evt_green_gcah2} illustrates that the new
algorithm requires only a few seconds to compute a highly accurate
approximation.

%
%

\section{Linear basis functions}

So far, we have only considered piecewise constant basis functions
in our experiments, since they make it particularly simple to
approximate the entries of the matrix $G$.
If the solution of the integral equation is smooth, it is generally
a good idea to employ basis functions of higher order to obtain
faster convergence.

One step up from piecewise constant basis functions are piecewise
linear functions, and we choose \emph{continuous} piecewise linear
functions, both in order to reduce the number of unknown variables
and to be able to work with integral operators that require an
$H^{1/2}$-conforming trial space.
The trial space is spanned by nodal basis functions $\varphi_i$:
$\varphi_i$ is continuous, piecewise linear on each triangle,
equal to one in the $i$-th vertex, and equal to zero in all
other vertices.

The support of $\varphi_i$ consists of all triangles that contain
the $i$-th vertex, therefore computing the entry $g_{ij}$ of the
matrix requires us to compute integrals on all pairs of triangles
$t\times s$ where $t$ belongs to the support of $\varphi_i$ and
$s$ to the support of $\varphi_j$:
\begin{equation*}
  g_{ij} = \sum_{t\subseteq\supp\varphi_i}
  \sum_{s\subseteq\supp\varphi_j}
  \int_t \varphi_i(x) \int_s g(x,y) \varphi_j(y) \,dy\,dx.
\end{equation*}
Due to this property, the computation of one entry of the matrix
with nodal basis functions can be \emph{significantly} more
computationally expensive than for a piecewise constant basis.

This problem can be somewhat mitigated by assembling the matrix
triangle pair by triangle pair:
we start with a zero matrix and loop over all pairs of triangles
$t\times s$.
For each pair, we evaluate the integrals for \emph{all} of the
triangles' vertices and add the results to the appropriate matrix
coefficients.
Although the final result is the same, we consider each
pair of triangles only once, and this allows us to re-use the
values of the kernel function $g$ in the quadrature points for
all combinations of basis functions.
Since the evaluation of the transformed kernel function is the
most computationally expensive part of the quadrature, this
approach can make the entire construction far more efficient.

Unfortunately, our compression scheme does not need \emph{all}
of the matrix entries, only entries for subsets
$\tilde\tau\times\tilde\sigma$ or $\hat\tau\times\hat\sigma$,
so looping over \emph{all} pairs of triangles would be a waste
of time.
Instead, we need an algorithm that determines only the required
triangles and loops over them.

We typically store a matrix $G|_{\hat\tau\times\hat\sigma}$ by enumerating
the row and column indices $\hat\tau=\{i_1,\ldots,i_n\}$,
$\hat\sigma=\{j_1,\ldots,j_m\}$ with $n=\#\hat\tau$,
$m=\#\hat\sigma$, and using a matrix in $\bbbr^{n\times m}$.
This means that it is not enough to find which triangles have to
participate in our computation, we also have to determine which
index numbers correspond to the triangles' vertices.

Given a standard representation of the mesh, it is quite simple
to determine for each index $i$ the set $T_i$ of triangles covering
the support of the basis function $\varphi_i$.
The challenge is to unify these sets for all basis functions
corresponding to a subset $\hat\tau$ of indices.
We use a variant of the well-known \emph{mergesort} algorithm
to handle this task:
for example, assume that we have triangles
\begin{align*}
  t_1 &= (1,2,3), &
  t_2 &= (2,3,5), &
  t_3 &= (4,1,3),\\
  t_4 &= (6,5,2), &
  t_5 &= (1,7,4), &
  t_6 &= (7,6,1)
\end{align*}
and are looking for the list of triangles for the vertex
set $\hat\tau=\{1,6,4\}$.
We have
\begin{align*}
  T_1 &= \{ t_1, t_3, t_5, t_6 \}, &
  T_6 &= \{ t_4, t_6 \}, &
  T_4 &= \{ t_3, t_5 \}.
\end{align*}
We write the triangles for each vertex into rows of a matrix,
where each row starts with the triangle, followed by three
entries for its three vertices that are equal to the local
index if this vertex is the current one or equal to the
special symbol $\bot$ if it is not:
\begin{equation*}
  \begin{array}{r|rrr}
    t_1 & 1 & \bot & \bot\\
    t_3 & \bot & 1 & \bot\\
    t_5 & 1 & \bot & \bot\\
    t_6 & \bot & \bot & 1\\
    \hline
    t_4 & 2 & \bot & \bot\\
    t_6 & \bot & 2 & \bot\\
    \hline
    t_3 & 3 & \bot & \bot\\
    t_5 & \bot & \bot & 3
  \end{array}
\end{equation*}
Now we apply the mergesort algorithm to sort the rows by the
first column.
If two rows have the same first column, i.e., if they correspond to
the same triangle, the rows are \emph{combined}:
if a column has an index in one row and $\bot$ in the other,
the combined row will have the index in this column.
If the rows have $\bot$ in the same column, the combined row will, too.
In our example, the result looks as follows:
\begin{equation*}
  \begin{array}{r|rrr}
    t_1 & 1 & \bot & \bot\\
    t_3 & 3 & 1 & \bot\\
    t_4 & 2 & \bot & \bot\\
    t_5 & 1 & \bot & 3\\
    t_6 & \bot & 2 & 1
  \end{array}
\end{equation*}
Each triangle appears in exactly one row, and each row
provides us with the local indices for all vertices of this
triangle.
The mergesort algorithm has a complexity of $\mathcal{O}(n \log n)$
if $n$ indices with $\mathcal{O}(1)$ triangles per index are used,
therefore the overhead for finding the triangles and the local indices
is low compared to the computational work for the quadrature itself.

\section{Curved triangles}

In our examples, linear basis function by themselves did reduce
the storage requirements, but did not lead to faster convergence
of the solution.
Since the reason appears to be that the polygonal approximation of
the smooth surface is insufficiently accurate, we consider replacing the
piecewise linear parametrizations of the triangles by piecewise
quadratic functions.
This leads to \emph{curved} triangles.

We implement these generalized triangles using the reference
triangle $\hat t := \{ x\in\bbbr^2\ :\ x_1,x_2\geq 0,\ x_1+x_2\leq 1 \}$
and quadratic parametrizations $\Phi_t,\Phi_s:\hat t\to\bbbr^3$
such that
\begin{align*}
  \int_t &\varphi_i(x) \int_s g(x,y) \varphi_j(y) \,dy \,dx\\
  &= \int_{\hat t} \gamma_t(\hat x) \varphi_i(\Phi_t(\hat x))
  \int_{\hat t} g(\Phi_t(\hat x),
  \Phi_s(\hat y))
  \gamma_s(\hat y) \varphi_j(\Phi_s(\hat y)) \,d\hat y \,d\hat x
\end{align*}
holds with the Gramians
\begin{align*}
  \gamma_t(\hat x)
  &= \left\|\frac{\partial \Phi_t}{\partial \hat x_1}(\hat x)
       \times \frac{\partial \Phi_t}{\partial \hat x_2}(\hat x)\right\|_2, &
  \gamma_s(\hat y)
  &= \left\|\frac{\partial \Phi_s}{\partial \hat y_1}(\hat y)
  \times \frac{\partial \Phi_s}{\partial \hat y_2}(\hat y)\right\|_2.
\end{align*}
For the basis functions, we choose \emph{mapped} nodal basis
functions, i.e., $\varphi_i\circ\Phi_t$ and $\varphi_j\circ\Phi_s$
are nodal linear basis functions on the reference triangle $\hat t$,
while $\varphi_i$ and $\varphi_j$ not necessarily linear themselves.

We can evaluate the double integral by using Sauter's quadrature
rule \cite{SA96,SASC11}, we only have to provide an efficient
way of evaluating the parametrization and the Gramian in the
quadrature points.
For the parametrization, we simply use quadratic interpolation
in the vertices and the midpoints of the edges.
For the Gramian, we observe that the outer normal vector
\begin{equation*}
  n_t(\hat x)
  = \frac{\partial \Phi_t}{\partial \hat x_1}(\hat x)
  \times \frac{\partial \Phi_t}{\partial \hat x_2}(\hat x)
\end{equation*}
is again a quadratic polynomial, so we can evaluate it also by
interpolation once we have computed its values in the vertices
and the midpoints.
Once we have $n_t(\hat x)$ at our disposal, computing
$\gamma_t(\hat x) = \|n_t(\hat x)\|_2$ is straightforward.
If we want to evaluate the double-layer potential operator and
need the \emph{unit} outer normal vector, we can obtain it by
simply dividing $n_t(\hat x)$ by $\gamma_t(\hat x)$.

%
%
\begin{figure}
  \begin{center}
  \includegraphics[width=0.8\textwidth]{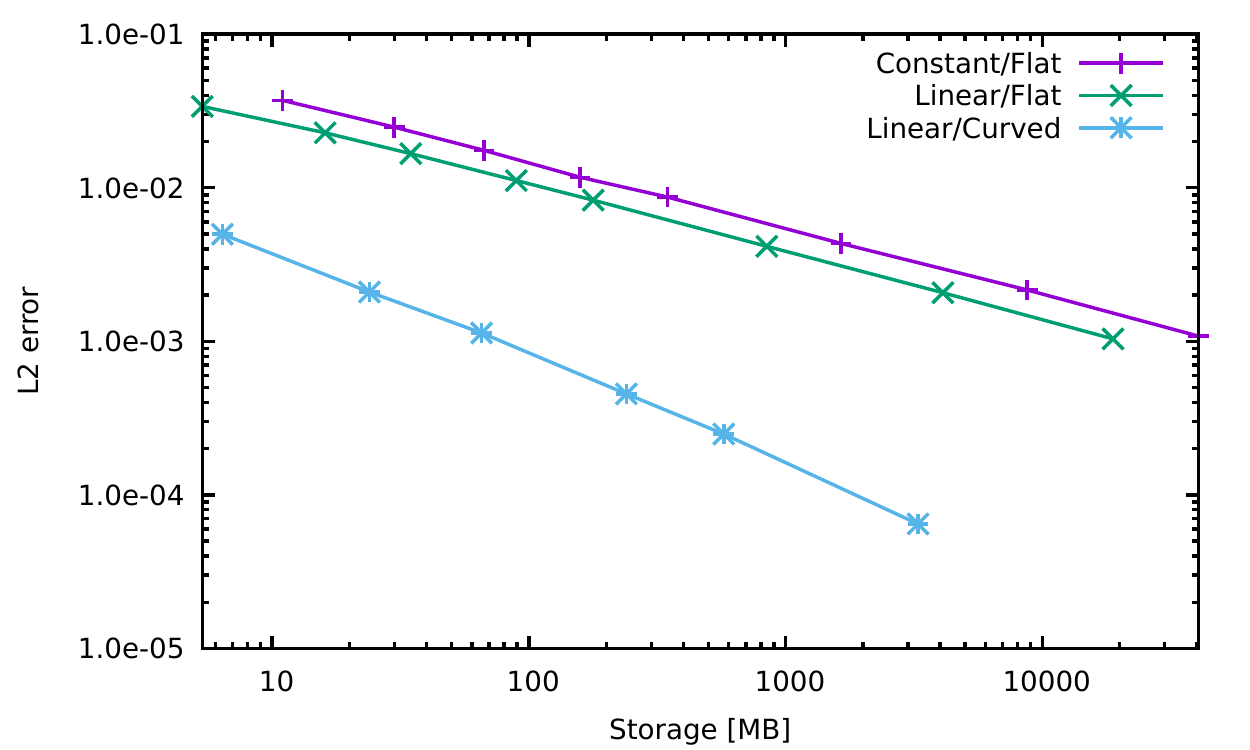}
  \caption{$L^2$ error compared to the storage requirements for
    constant and linear basis functions as well as plane and
    curved triangles}
  \label{fi:evm_cc_ll}
  \end{center}
\end{figure}

Figure~\ref{fi:evm_cc_ll} shows the $L^2$ error of the solution compared
to the storage requirements for constant basis functions, linear basis
functions on plane triangles, and linear basis functions on curved triangles.
We can see that constant and linear basis functions on plane
triangles converge at approximately the same rate, while curved
triangles lead to a significantly improved rate of convergence
and pronouncedly smaller errors for identical problem dimensions.

%
%
\begin{figure}
  \begin{center}
  \includegraphics[width=0.8\textwidth]{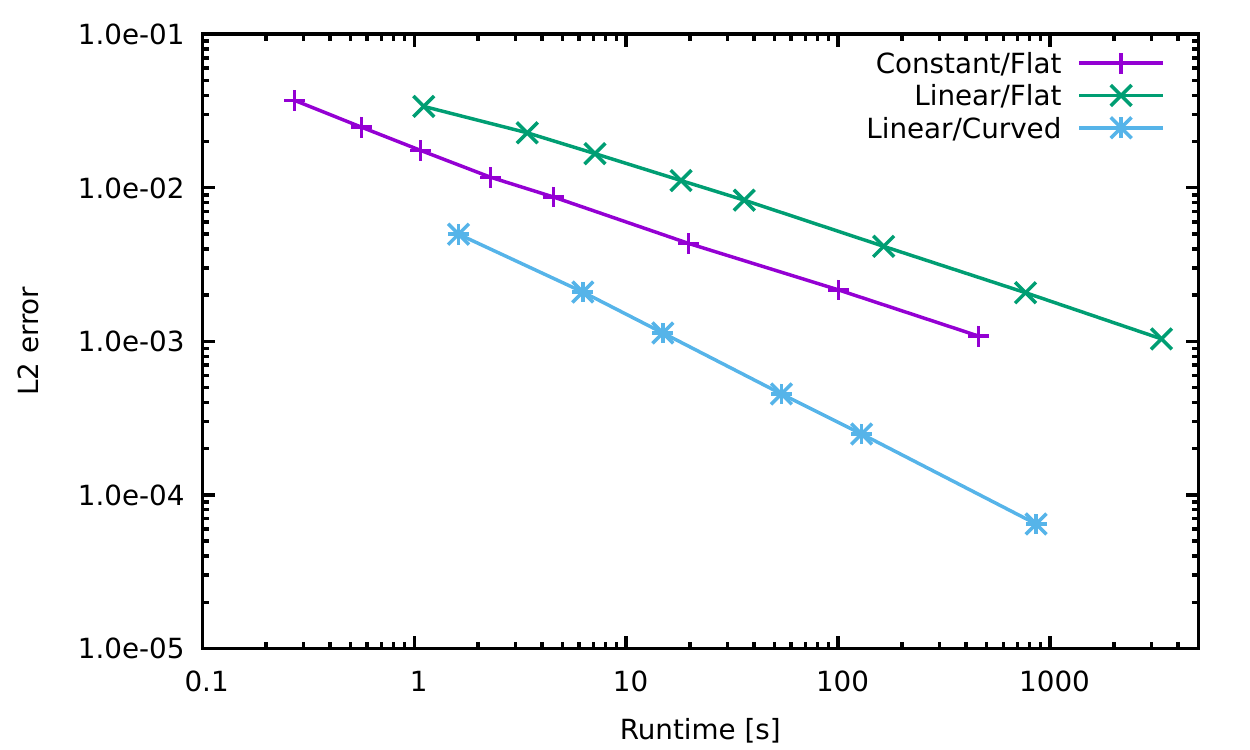}
  \caption{$L^2$ error compared to the setup time for
    constant and linear basis functions as well as plane and
    curved triangles}
  \label{fi:evt_cc_ll}
  \end{center}
\end{figure}

Figure~\ref{fi:evt_cc_ll} compares the $L^2$ error to the setup times
for the three cases.
As is to be expected, linear basis functions require far more
time than constant basis functions, and curved triangles again
take more time than plane ones.
But we can also see that the combination of linear basis
functions with curved triangles provides us with significantly
lower $L^2$ errors, making the most sophisticated approach also the
most efficient of the three.

If the surface and the solution are sufficiently smooth, we can choose
\emph{collocation} instead of Galerkin's method for the
discretization, i.e., define the matrix entries via
\begin{equation*}
  g_{ij} = \int_{\partial\Omega} g(x_i,y) \varphi_j(y) \,dy,
\end{equation*}
where $x_i$ is a vertex of the mesh and $\varphi_j$ is a nodal
basis function.
Since collocation requires only a single instead of a double
integral, the number of quadrature points is significantly
smaller.
Since the singularity is fixed at $y=x_i$, a simple Duffy
transformation is sufficient to regularize the integral, and this
makes the implementation quite straightforward.

%
%
\begin{figure}
  \begin{center}
  \includegraphics[width=0.8\textwidth]{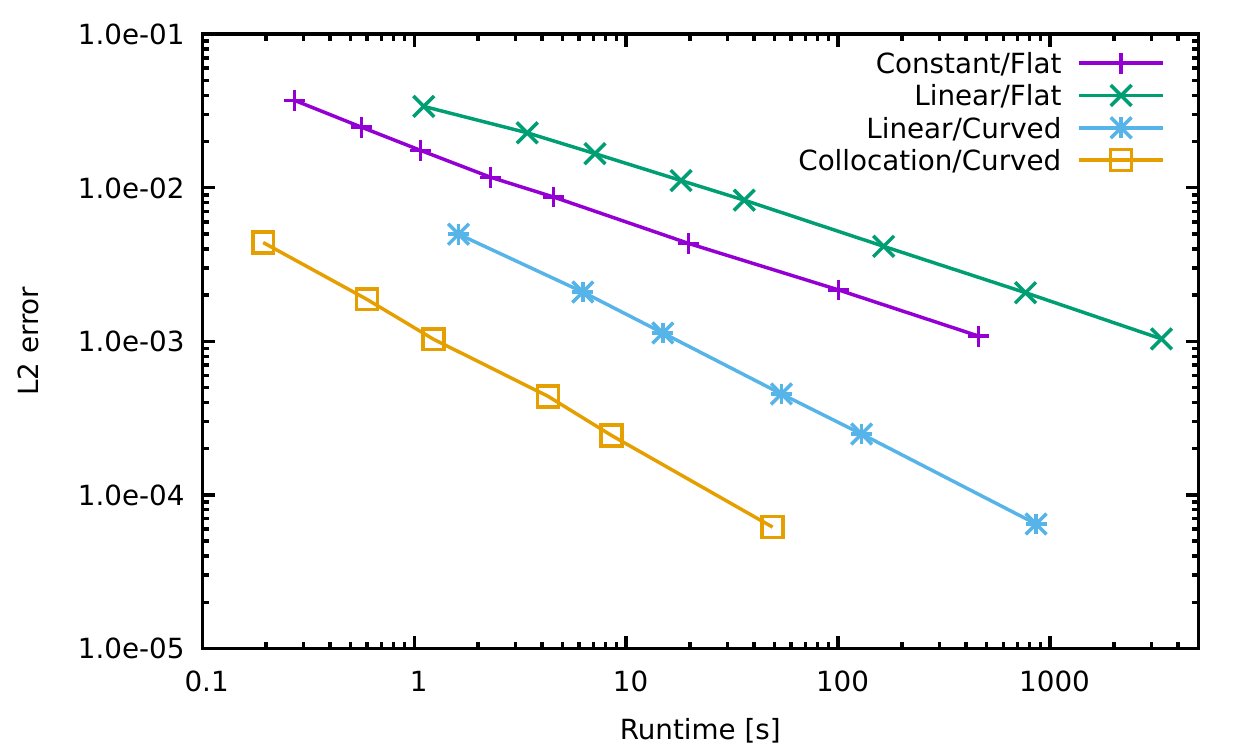}
  \caption{$L^2$ error compared to the setup time for
    constant and linear basis functions as well as plane and
    curved triangles}
  \label{fi:evt_coll}
  \end{center}
\end{figure}

Figure~\ref{fi:evt_coll} illustrates that collocation reduces
the setup time by a factor of approximately ten compared to
the Galerkin discretization without significantly changing the
quality of the approximated solution.

\section{GPU implementation}

Modern computers are frequently equipped with powerful graphics
processors that are (reasonably) pro\-grammable and can therefore
help with certain computational tasks.
These processors are frequently called \emph{general-purpose
  graphics processing units} (GPGPUs or short GPUs) and differ
substantially from standard processors (CPUs).
In order to use GPUs to accelerate our algorithm, we have to
take the architectural differences into account.

A first important difference is the way CPUs and GPUs handle
data:
high-end GPUs typically are connected to dedicated high-bandwidth
memory.
While a current CPU may reach a memory bandwidth of $60$ GBytes/s,
modern GPUs provide up to $550$ GBytes/s.
It has to be pointed out that the higher bandwidth comes at a price:
while even desktop CPUs can access $64$ GBytes of RAM, with server
CPUs accessing up to $1024$ GBytes, current GPUs are limited to
$24$ GBytes of memory.
In order to deal with large data sets, we have to move data between
graphics memory and main memory, and these transfers are fairly slow.

The most important difference is the number of arithmetic units:
while a 28-core CPU with 512-bit vector registers can
perform $28\times 16=448$ double-precision floating-point operations
per clock, high-end GPUS offer currently up to $4\,608$ arithmetic
units that can work in parallel.
Even taking differences in clock speeds into account, the theoretical
computing power of GPUs is significantly larger than that of CPUs.

In order to control the large number of arithmetic units efficiently,
GPUs restrict the ways these units can work.
Current architectures follow what is known as the \emph{single instruction,
multiple threads} (SIMT) model:
the computation is split into \emph{threads}, frequently hundreds of
thousands or millions, each with its own instruction pointer and
local variables.

In order to keep the management of the threads simple, a fixed number
of threads is bundled into a \emph{warp}, e.g., $32$ or $64$ threads,
depending on the architecture.

The GPU consists of multiple \emph{multiprocessors} that can execute
the instructions required by a warp.
Each multiprocessor is assigned a certain number of warps.
In each cycle, one of these warps and one of its instructions is
chosen for execution.
If the instruction pointer of a thread indicates the chosen instruction,
it is executed, otherwise the thread remains idle during the current
cycle.

This is a key difference between GPUs and CPUs:
all threads running on a CPU are independent and can execute any
instruction per cycle, while all threads in the same warp on
a GPU have to execute the same instruction or do nothing.

If the control flow of the threads within one warp diverges,
i.e., if all of the threads have to execute different instructions,
only one of the instruction can be executed per cycle, allowing
only one of the threads to advance.
Obviously, having $63$ of $64$ arithmetic units idle for an
extended period of time is not the best use of the available
hardware.

\section{\texorpdfstring{GCA-$\mathcal{H}^2$ for GPUs}{GCA-H2 for GPUs}}

Let us now consider how to adapt our algorithm for execution on
GPUs.

The computational work is dominated by three tasks:
\begin{itemize}
  \item the construction of the leaf and transfer matrices $V_\tau$
    and $E_\tau$ and the index sets $\tilde\tau$ by Green quadrature
    and cross approximation,
  \item the computation of the coupling matrices
    $G|_{\tilde\tau\times\tilde\sigma}$ for admissible blocks, and
  \item the computation of $G|_{\hat\tau\times\hat\sigma}$ for the
    remaining inadmissible blocks.
\end{itemize}
Although the entries of $A_{\tau\sigma}$ or $\widehat{A}_{\tau\sigma}$
involve no control-flow divergence and should therefore be well-suited
for SIMT architectures, the highly adaptive nature of the cross
approximation leads us to leave the first part of the algorithm
to the CPU, where parallelization and vectorization can be employed
to take full advantage of the available resources.

Once the sets $\tilde\tau$ and $\tilde\sigma$ are known, the computation
of the entries of the matrices $G|_{\tilde\tau\times\tilde\sigma}$ and
$G|_{\hat\tau\times\hat\sigma}$ requires no adaptivity whatsoever, so the
second and third part of our algorithm can be expected to be ideally
suited for GPUs.

Setting up the GPU to run a number of threads involves a certain
amount of communication and management operations, therefore we should
make sure that the number of threads is sufficiently high in order to
minimize organizational overhead.
Since our algorithm only works with small matrices that would not
allow us to reach an adequate number of threads, we switch to an
asynchronous execution model:
instead of computing the entries of a matrix the moment it is
encountered by our algorithm, the corresponding task is added to
a list for later handling.
Only once the list has grown enough to keep a sufficiently large
number of threads busy, it is transferred to the GPU for execution.

This approach also allows us to handle different cases appearing
in the numerical quadrature:
we use Sauter's quadrature rule \cite{SA96,SASC11} to integrate
the singular kernel function on pairs of triangles.
Sauter's algorithm requires different quadrature points (and even
different numbers of quadrature points) depending on whether the
triangles are identical, share an edge, a vertex, or are disjoint.
By simply using one list for each of the four cases, we can ensure
that all threads execute almost exactly the same sequence of
instructions and that control-flow divergence is kept down to
a minimum.

Since communication between the CPU and the GPU is slow, we should
try to keep the amount of data that has to be transferred as small
as possible.
In our implementation, the geometrical information of the triangles
is kept permanently in graphics memory, so that we only have to
transfer the numbers of the triangles $t$ and $s$ in order to describe
an integral that has to be computed.

Another important step in reducing the impact of the communication
between CPU and GPU is to ``hide'' the communication behind
computation:
modern graphics cards can perform computations and memory transfers
concurrently, and we use this feature in order to use the time
spent by the arithmetic units on one list to transfer the results
of the previous list back to main memory and the input of the next
list to graphics memory.

Finally, we employ multiple threads on the CPU to fill multiple
lists concurrently in order to ensure that both the memory
management and the arithmetic units of the GPU are kept busy.

%
%
\begin{figure}
  \begin{center}
  \includegraphics[width=0.8\textwidth]{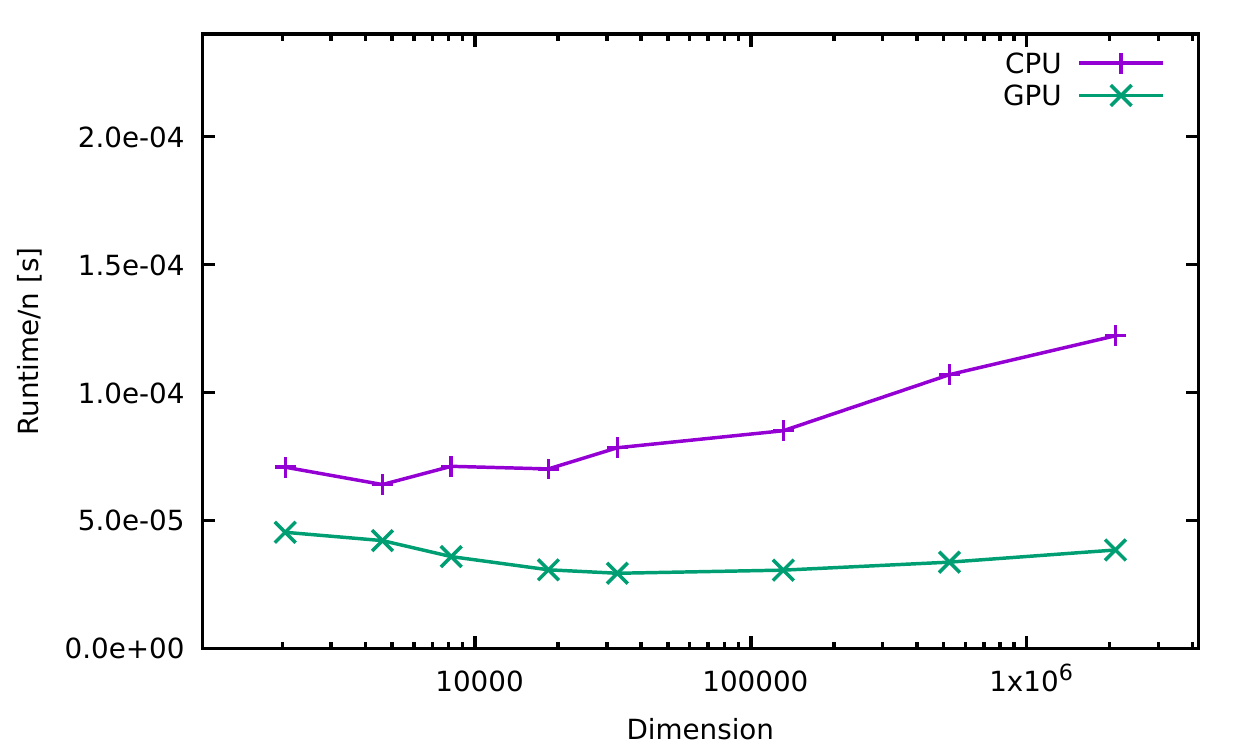}
  \caption{Runtime of CPU and GPU setup of the GCA-$\mathcal{H}^2$-matrix}
  \label{fi:cpu_gpu_flat}
  \end{center}
\end{figure}

Figure~\ref{fi:cpu_gpu_flat} shows the runtime in seconds per degree
of freedom for setting up the GCA-$\mathcal{H}^2$ matrix for different
meshes approximating the unit sphere.
For the CPU, we use an Intel Core i7-7820 with 8 cores and AVX 512
running at a base frequency of $3.6$ GHz, providing a peak performance
of $460$ GFlops/s at double precision.
For the GPU, we have used an AMD Vega 64 card running at 1.25 GHz
with 8 GBytes of HBM2 memory and $4\,096$ arithmetic units providing
a peak performance of 791 GFlops/s at double precision (and
\emph{considerably} more for single precision).
It should be pointed out that this GPU is designed for the consumer
market, so its double-precision performance is quite low.
GPUs designed for computation like the NVIDIA Tesla P100 or V100
should provide around $5\,000$ and $7\,000$ GFlops/s at double
precision, respectively.

Our figure suggests that both the CPU and the GPU implementation
have $\mathcal{O}(n \log n)$ complexity:
we use a logarithmic scale for the dimension $n$ and a linear
scale for the work per degree of freedom, and the figure shows the
latter as, essentially, a line.
We can also see that the slope of this line for the CPU implementation
is significantly steeper than for the GPU implementation.
This suggests that the GPU implementation may become increasingly
more attractive for larger meshes.

\bibliographystyle{plain}
\bibliography{hmatrix}

\begin{thebibliography}{10}

\bibitem{BE00a}
M.~Bebendorf.
\newblock Approximation of boundary element matrices.
\newblock {\em {N}umer. {M}ath.}, 86(4):565--589, 2000.

\bibitem{BERJ01}
M.~Bebendorf and S.~Rjasanow.
\newblock {A}daptive {L}ow-{R}ank {A}pproximation of {C}ollocation {M}atrices.
\newblock {\em Computing}, 70(1):1--24, 2003.

\bibitem{BO10}
S.~{B\"orm}.
\newblock {\em Efficient Numerical Methods for Non-local Operators: {${\mathcal
  H}^2$}-Matrix Compression, Algorithms and Analysis}, volume~14 of {\em EMS
  Tracts in Mathematics}.
\newblock EMS, 2010.

\bibitem{BOCH14}
S.~{B\"orm} and S.~Christophersen.
\newblock Approximation of integral operators by {Green} quadrature and nested
  cross approximation.
\newblock {\em Numer. Math.}, 133(3):409--442, 2016.

\bibitem{BOGR02}
S.~{B\"orm} and L.~Grasedyck.
\newblock Low-rank approximation of integral operators by interpolation.
\newblock {\em Computing}, 72:325--332, 2004.

\bibitem{BOGR04}
S.~{B\"orm} and L.~Grasedyck.
\newblock Hybrid cross approximation of integral operators.
\newblock {\em Numer. Math.}, 101:221--249, 2005.

\bibitem{BOHA02}
S.~{B\"orm} and W.~Hackbusch.
\newblock Data-sparse approximation by adaptive {${\mathcal{H}}^2$}-matrices.
\newblock {\em Computing}, 69:1--35, 2002.

\bibitem{BOLOME02}
S.~{B\"orm}, M.~L{\"o}hndorf, and J.~M. Melenk.
\newblock Approximation of integral operators by variable-order interpolation.
\newblock {\em Numer. Math.}, 99(4):605--643, 2005.

\bibitem{CHIP94}
S.~Chandrasekaran and I.~C.~F. Ipsen.
\newblock On rank-revealing factorisations.
\newblock {\em SIAM J. Matrix Anal. Appl.}, 15(2):592--622, 1994.

\bibitem{GIRO02}
Z.~Gimbutas and V.~Rokhlin.
\newblock A generalized fast multipole method for nonoscillatory kernels.
\newblock {\em SIAM J. Sci. Comput.}, 24(3):796--817, 2002.

\bibitem{GRRO97}
L.~Greengard and V.~Rokhlin.
\newblock A new version of the fast multipole method for the {L}aplace equation
  in three dimensions.
\newblock In {\em Acta Numerica 1997}, pages 229--269. Cambridge University
  Press, 1997.

\bibitem{HAKHSA00}
W.~Hackbusch, B.~N. Khoromskij, and S.~A. Sauter.
\newblock On $\mathcal{H}^2$-matrices.
\newblock In H.~Bungartz, R.~Hoppe, and C.~Zenger, editors, {\em Lectures on
  Applied Mathematics}, pages 9--29. Springer-Verlag, Berlin, 2000.

\bibitem{HANO89}
W.~Hackbusch and Z.~P. Nowak.
\newblock On the fast matrix multiplication in the boundary element method by
  panel clustering.
\newblock {\em Numer. Math.}, 54(4):463--491, 1989.

\bibitem{RO85}
V.~Rokhlin.
\newblock Rapid solution of integral equations of classical potential theory.
\newblock {\em J. Comp. Phys.}, 60:187--207, 1985.

\bibitem{SA96}
S.~A. Sauter.
\newblock Cubature techniques for 3-d {G}alerkin {BEM}.
\newblock In W.~Hackbusch and G.~Wittum, editors, {\em Boundary Elements:
  Implementation and Analysis of Advanced Algorithms}, pages 29--44.
  Vieweg-Verlag, 1996.

\bibitem{SASC11}
S.~A. Sauter and C.~Schwab.
\newblock {\em Boundary Element Methods}.
\newblock Springer, 2011.

\bibitem{TY96}
E.~E. Tyrtyshnikov.
\newblock Mosaic-skeleton approximation.
\newblock {\em {C}alcolo}, 33:47--57, 1996.

\bibitem{BIYIZO04}
L.~Ying, G.~Biros, and D.~Zorin.
\newblock A kernel-independent adaptive fast multipole algorithm in two and
  three dimensions.
\newblock {\em J. Comp. Phys.}, 196(2):591--626, 2004.

\end{thebibliography}

\end{document}